\documentclass[11pt,a4paper]{article}
\usepackage{amsmath,amstext,amssymb,amscd}
\newcommand{\Xcomment}[1]{}

\makeatletter \@addtoreset{equation}{section} \makeatother

{\hfill$\qed$\medskip}

\def\qed{ \ \vrule width.1cm height.3cm depth0cm}

 \makeatletter
\renewcommand{\section}{\@startsection{section}{1}{0pt}%
{-3.5ex plus -1ex minus -.2ex}{2.3ex plus .2ex}%
{\normalfont\Large}}
 \makeatother

 \makeatletter
\renewcommand{\subsection}{\@startsection{subsection}{2}{0pt}%
{-3.0ex plus -1ex minus -.2ex}{1.5ex plus .2ex}%
{\normalfont\normalsize\bf}}
 \makeatother

\input xy
\xyoption{all}

\begin{document}
\title{
Cluster monomials in $\mathbb C[GL_n/N]$, a simplicial fan in the cone of semi-standard Young tableaux, and the Lusztig basis}

\author{G.A.Koshevoy\footnote{\noindent CEMI and Poncelet laboratoty (IMU and CNRS (UMI 2615)),
email: koshevoy@cemi.rssi.ru} }

\date{}
\maketitle
\begin{abstract}
We study the cluster monomials and cluster complex in $\mathbb C[GL_n/N]$. For we consider the {\em tableau basis} in $\mathbb C[GL_n/N]$. Namely, an element $\Delta_T$ of the tableau basis labeled by a semistandard Young tableau $T$ is the product of the flag minors corresponding to columns of $T$.
Our main results state: (i) cluster monomials in $\mathbb C[GL_n/N]$ can be labeled by semistandard Young tableaux such that any cluster monomial has the form $\Delta_T+$ lexicographically smaller terms;
(ii) such labeling distinguish the cluster monomials; (iii) for any seed of the cluster algebra on $\mathbb C[GL_n/N]$,
we define a cone in $\mathbf D(n)$  ($\mathbf D(n)$ is the cone of semi-standard Young tableaux,
$\mathbf D(n)$ is linear isomorphic to the Gelfand-Tseitlin cone)
generated by tableaux which label the cluster variables of the seed, and these cones form a simplicial fan in
$\mathbf D(n)$.
\end{abstract}

\section{Introduction}
One of the main motivation of S.Fomin and A.Zelevinsky
for introducing cluster algebras \cite{FZ1} was the desire to
provide a combinatorial framework to understand the structure of 'dual
canonical bases' in (homogeneous) coordinate rings of various algebraic varieties related to semisimple groups.
Several such varieties (Grassmann varieties and double Bruhat cells)
carry a cluster algebra structure and certain
special functions on that spaces (such as Pl\"ucker coordinates, generalized minors) correspond to distinguished elements
called {\em cluster variables}.

For the general linear group $G=GL_n(\mathbb C)$, and the subgroup $N\subset G$ of unipotent
upper-triangular matrices,  the coordinate ring of the {\em  base affine space} for $GL_n(\mathbb C)$,
$\mathbb C[GL_n/N]$, is
the set of regular functions on $GL_n(\mathbb C)$ which are invariant under the action of $N$
by the right multiplication.
According to classical invariant theory, $\mathbb C[GL_n/N]$ is
generated by the {\em flag minors}. (A flag minor $\Delta_I$ of a matrix $x=(x_{ij})\in GL_n(\mathbb C)$ is a minor
occupying in rows in $I$ and the first $|I|$
columns.)  The ring $\mathbb C[GL_n/N]$ is one of prototypical examples of a cluster algebra (\cite{Fomin10})).
For the cases $n=3,4,5$, the cluster monomials form a basis in  $\mathbb C[GL_n/N]$ as in the vector space,  but
for $n\ge 6$ the cluster monomials  fail to be a basis.

We study the cluster monomials and cluster complex in $\mathbb C[GL_n/N]$.
For we consider the {\em tableau basis} in $\mathbb C[GL_n/N]$
(see, for example \cite{DKR,Fulton}).
The elements of the tableau basis are labeled by the semistandard Young tableaux.
Namely, an element $\Delta_T$
of the tableau basis labeled a semistandard Young tableau $T$ is the
product of the flag minors
corresponding to columns of $T$ (for details see Section 3). Because of the   D\'esam\'enien-Kung-Rota algorithm (\cite{DKR}), the set of such monomials in flag minors is a basis in $\mathbb C[GL_n/N]$.

Consider the lexicographical order $\succ$ on the variables $x_{ij}$, $i, j\in  [n]$,
($i$ indicates the row number and $j$ the column number) $x_{ij}\succ x_{kl}$ if
$i<k$ or $i=k$ and $j<l$.  This order defines a total order on the monomials of $\mathbb C[x_{ij}]$.
In particular, we get a total order on the Young tableaux through the total ordering of $\Delta_T$'s.

The semi-standard Young tableaux filled in the alphabet $[n]$
can be viewed as integer point of the cone of $\mathbf D$-tight arrays
$\mathbf D(n)\subset\mathbb R^{n(n+1)/2}$ (see \cite{array}). The cone $\mathbf D(n)$ is linear isomorphic
to well-known cone of the Gelfand-Tseitlin patterns. For us it is convenient to
work with the cone of semi-standard Young tableaux.

Our main results (Theorem M1, M2 and M3) state:
(i) cluster monomials in $\mathbb C[GL_n/N]$ can be labeled by semistandard Young tableaux
such that any cluster monomial has the form $\Delta_T+$ lexicographically smaller terms;
(ii) such labeling distinguish the cluster monomials; (iii) for any seed of the cluster algebra on
$\mathbb C[GL_n/N]$, we define a cone in $\mathbf D(n)$ generated by tableaux which label the cluster
variables of the seed; these cones form a simplicial fan in $\mathbf D(n)$.

In the case the cluster algebra $\mathbb C[GL_n/N]$ is of finite type, that is the case for $n=3,4,5$, the union of
the cones of the fan coincides with the cone $\mathbf D(n)$, and thus the cluster
monomials form a basis in $\mathbb C[GL_n/N]$.

For $n\ge 6$, the cluster algebra of $\mathbb C[GL_n/N]$ is of infinite type, and the union
of the cones is only a part of the cone $\mathbf D(n)$.

We conjecture that the cluster monomials correspond to the {\em real} elements of the Lusztig basis.
Namely,
the specialization of the dual canonical basis at $q=1$ is the {\em Lusztig basis}. Recall, that the dual canonical basis is a basis of the quantum deformation of $\mathbb  C[GL_n/N]$.
This is a distinguished basis which nicely behaves with respect to implementation of $\mathbb C[GL_n/N]$ as a representation of $GL_n$ (namely, $\mathbb C[GL_n/N]$ is the direct sum of the irreducible representations of $GL_n$ each taken with multiplicity one).

The Lusztig basis is a basis in $\mathbb C[GL_n/N]$ labeled by the integer point of the cone $\mathbf D(n)$.
An element of Lusztig basis labeled by an integer point of $\mathbf D(n)$ has the form
$\Delta_T+$ lexicographically smaller terms, where $T$ is a semi-standard Young tableau corresponding to the point
in $\mathbf D(n)$.

The Lusztig basis and the tableau basis are different.
For example, for $GL_3$,  the Lusztig basis (discovered in  1985 by I.Gelfand and A.Zelevinsky
(\cite{GZ})) is the collection of monomials in $\Delta_{1}$, $\Delta_{2}$,
$\Delta_{3}$, $\Delta_{12}$, $\Delta_{23}$, $\Delta_{13}$, and $\Delta_{123}$ which do not contain the product
$\Delta_2\Delta_{13}$, while                                                                                              the tableaux basis is the collection of monomials in $\Delta_{1}$, $\Delta_{2}$, $\Delta_{3}$, $\Delta_{12}$,
$\Delta_{23}$, $\Delta_{13}$, and $\Delta_{123}$ which do not contain the monomial $\Delta_1\Delta_{23}$.

An element of the Lusztig basis is {\em  real} if
its square belongs to the Lusztig basis.  We call the semi-standard Young tableau {\em real}
if the corresponding element of the Lusztig basis is real.
Our conjecture is that cluster monomials are labeled by real semi-standard Young tableaux.
\medskip

{\em Acknowledgments}. I thank Vladimir Danilov, Sergey Fomin, Alexander Karzanov, and Jan Schr\"oer
for useful discussions. A part of this work was made in the Max-Planck
Institute f\"ur Mathematics,  Hausdorff Institute for Mathematics (Bonn) and IHES (Bures-sur-Yvette) and I
thanks these institutes
for hospitality and financial support.

\section{Cluster skew-symmetric algebras}
Since we are interested in the cluster algebra structure on $\mathbb C[GL_n/N]$,
we remind necessary definitions for so-called skew-symmetric
cluster algebras.

Let $G=(V(G),E(G))$ be a {\em quiver} (a directed multigraph) in which the vertex
set $V(G)$ is partitioned into two subsets: a set $V_1$ of
\emph{frozen} vertices, and a set $V_2$ of \emph{mutable}
vertices. The (integer) edge multiplicity function is regarded as
being skew-symmetric: if vertices $u,v$ are connected by $\alpha$
edges going from $u$ to $v$ (which are members of $E(G)$), we
simultaneously think of these vertices as being connected by
$-\alpha$ edges going from $v$ to $u$. To each vertex $v$ of $G$
one associates a {\em cluster  variable} $x_v$ so that $\{x_v\colon v\in V(G)\}$
is a transcendence basis of a field of rational functions $\mathcal F$. Such a
pair consisting of a quiver and a transcendence basis indexed by
its vertices is said to be a {\em cluster seed\/}. Monomials in $x_v$, $v\in V(G)$ are called {\em cluster monomials}.
The cluster seed can be mutated at any mutable vertex to produce a new
cluster variable and a new seed. Applying mutations in all possible situations produce the set of cluster variables and these variables form a
skew-symmetric cluster algebra~\cite{FZ1}.

The quiver and variables are modified by applying the following
operations called cluster mutations. A \emph{cluster mutation}
$\mu_v$ applied at a mutable vertex $v\in V_2$ changes one
variable, namely, $x_v$, and modifies the quiver $G$ by the following rules.
For a vertex $v$, denote $In(v):=\{v'\in V(G)\,:\, (v',v)\in
E(G)\}$ and $Out(v):=\{v''\in V(G)\,:\, (v,v'')\in E(G)\}$.

The quiver $\mu_v(G)$ has the same vertex set as $G$,
$V(\mu_v(G))=V(G)$, partitioned into frozen and mutable vertices
in the same way as before. The edges $E(\mu_v(G))$ are obtained
from edges $E(G)$ by the following rule:

(i) the edges in $E(\mu_v(G))$ incident to the vertex $v$ are
exactly the edges in $E(G)$ incident to $v$ but taken with the
reverse direction;

(ii) for each pair $v'\in In(v)$ and $v''\in Out(v)$, form the
edge $(v',v'')$ in $E(\mu_v(G))$ whose multiplicity is defined to
be $\gamma-\alpha\cdot\beta$, where $\alpha=w(v',v)\ge 1$ is multiplicity
of the edge $(v',v)$ in $E(G)$, $\beta=w(v,v'')\ge 1$ is that for
$(v,v'')$, and $\gamma\in\mathbb Z$ is that for $(v',v'')$;

(iii) the other edges of $\mu_v(G)$ are those of $G$ that neither
are incident to $v$ nor connect pairs $v',v''$ as in~(ii).\medskip

For $u\neq v$, we put $\mu_v(x_u):=x_u$ and define
$\mu_v(x_v)=x^{new}_v$ by the following rule:
  $$
x^{new}_v\cdot x_v=\prod_{v'\in In(v)}x_{v'}^{w(v',v)}+\prod_{v''\in
Out(v)}x_{v''}^{w(v,v'')}.
  $$
This gives the new seed: the quiver $\mu_v(G)$ and variables $\mu_v(x_u)$,
~$u\in V(\mu_v(G))=V(G)$.

Obviously, there holds $\mu_v^2=id$.

The variables $x_v$, $v\in V_1$ do not change and are called {\em coefficients}.  Let $\mathcal X$ denotes the set of
cluster variables. Then the cluster algebra is $k[x_v\,|\,v\in V_1]$-subalgebra of $\mathcal F$ generated by $\mathcal X$.

We need the following important results on cluster algebras. Namely, the Laurent phenomenon (\cite{FZ1}) and recently proven Positivity conjecture for
skew-symmetric cluster algebras \cite{pos}.
\medskip

{\bf Theorem CL}. (\cite{FZ1}) For any initial seed $(G;x_v,\, v\in V(G))$,
any cluster variable $x\in\mathcal X$ is a Laurent polynomial in variables $x_v, v\in V$.\medskip

{\bf Theorem CP}. (\cite{pos}) For a skew-symmetric cluster algebra, the coefficients in such Laurent polynomials are positive.\medskip

\section{Cluster monomials in $\mathbb C[GL_n/N]$ and a fan in the cone of semi-standard Young tableaux}

For given $n$, let us consider cluster algebras $\mathcal A(n)$ with the initial seed specified by
the quiver $Q(n)$ being a
triangular grid of size $n$ with cyclically oriented triangles (we depicted $Q(5)$ below).

We are interested in such  a
cluster algebra, because the quiver $Q(n)$ corresponds to a pseudoline
arrangement for the reduced decomposition of the longest permutation  \\
$w_0=s_1s_2\ldots s_{n-1}s_1\ldots s_{n-2}\ldots s_1s_2s_1$ (see \cite{Fomin10} and  \cite{BFZup}).
\bigskip

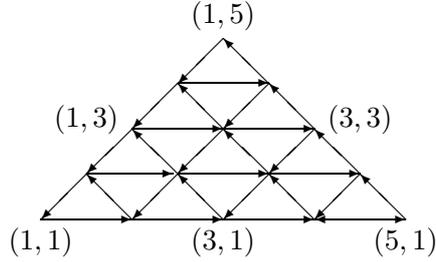
 \begin{center}
\unitlength=0.6mm \special{em:linewidth 0.4pt}
\linethickness{0.4pt}
\begin{picture}(100.00,60.00)
\put(30.00,30.00){\vector(-1,-1){10.00}}
\put(40.00,40.00){\vector(-1,-1){10.00}}
\put(50.00,50.00){\vector(-1,-1){10.00}}
\put(60.00,60.00){\vector(-1,-1){10.00}}
\put(70.00,50.00){\vector(-1,1){10.00}}
\put(80.00,40.00){\vector(-1,1){10.00}}
\put(90.00,30.00){\vector(-1,1){10.00}}
\put(100.00,20.00){\vector(-1,1){10.00}}
\put(60.00,40.00){\vector(-1,1){10.00}}
\put(70.00,30.00){\vector(-1,1){10.00}}
\put(80.00,20.00){\vector(-1,1){10.00}}
\put(50.00,30.00){\vector(-1,1){10.00}}
\put(60.00,20.00){\vector(-1,1){10.00}}
\put(40.00,20.00){\vector(-1,1){10.00}}
\put(70.00,50.00){\vector(-1,-1){10.00}}
 \put(60.00,40.00){\vector(-1,-1){10.00}}
\put(50.00,30.00){\vector(-1,-1){10.00}}
\put(80.00,40.00){\vector(-1,-1){10.00}}
\put(70.00,30.00){\vector(-1,-1){10.00}}
\put(90.00,30.00){\vector(-1,-1){10.00}}
\put(100.00,20.00){\vector(-1,0){20.00}}
\put(80.00,20.00){\vector(1,0){20.00}}
\put(60.00,20.00){\vector(1,0){20.00}}
\put(40.00,20.00){\vector(1,0){20.00}}
\put(20.00,20.00){\vector(1,0){20.00}}
\put(70.00,30.00){\vector(1,0){20.00}}
\put(50.00,30.00){\vector(1,0){20.00}}
\put(30.00,30.00){\vector(1,0){19.00}}
\put(60.00,40.00){\vector(1,0){20.00}}
\put(40.00,40.00){\vector(1,0){20.00}}
\put(50.00,50.00){\vector(1,0){20.00}}
\put(20.00,15.00){\makebox(0,0)[cc]{$(1,1)$}}
\put(100.00,15.00){\makebox(0,0)[cc]{$(5,1)$}}
\put(60.00,65.00){\makebox(0,0)[cc]{$(1,5)$}}
\put(60.00,15.00){\makebox(0,0)[cc]{$(3,1)$}}
\put(90.00,42.00){\makebox(0,0)[cc]{$(3,3)$}}
\put(30.00,42.00){\makebox(0,0)[cc]{$(1,3)$}}
\end{picture}
 \end{center}

For us it is convenient to label the vertices of $Q(n)$ by  $(i,j)$, $i+j\le n+1$, $i$, $j\ge 1$,
such that
the vertex labeled by $(1,1)$ is at the left corner of the triangular  grid, and $(i,j)$ labels
the end point of any path in $Q$ which has $i-1$ horizontal edges directed to East and $j-1$
slope  edges directed North-East (note, that
the north-east direction is the opposite to the direction of the corresponding edges in $Q(n)$).
(We depicted the vertices (1,1), (3,1), (1,3), (1,5), (5,1) and (3,3) in the above picture.) The vertices at the left and right sides  of the triangle are frozen, that is the set $V_1$ of
the frozen vertices consists of vertices labeled by $(i,1)$, $i=1,\ldots ,n$, and $(i,j)$ with $i+j=n+1$,
and the set $V_2$
of mutable vertices is the vertices labeled by $\{(i,j)$, either $i\neq 1$ or $i+j<n+1\}$.

The cluster structure on $\mathbb C[GL_n/N]$ is a specialization of the cluster algebra  $\mathcal A(n)$.
Namely, we view $n(n+1)/2$  interval flag minors $\Delta_{\{i,i+1, \ldots, i+j-1\}}$
labeling the vertices $(i,j)$, $i+j\le n+1$, $i, j\ge 1$,
as formal indeterminates. There are $(n-2)(n-1)/2$ possible mutations out of the seed $Q(n)$,
we use the quiver to write  the corresponding exchange relations:
\[
\Delta_{\{i,\ldots, i+j-1\}}\Omega_{ij}= \Delta_{\{i,\ldots, i+j-2\}}\Delta_{\{i+1,\ldots, i+j\}}\Delta_{\{i-1,\ldots, i+j-1\}}+
\Delta_{\{i,\ldots, i+j\}}\Delta_{\{i+1,\ldots, i+j-1\}}\Delta_{\{i-1,\ldots, i+j-2\}}, j>1,
\]
and, for $j=1$, we have
\[
\Delta_i\Omega_i= \Delta_{i-1}\Delta_{i,i+1}+\Delta_{i-1,i}\Delta_{i+1}, i\ge 2.
\]
Using the Pl\"ucker relations, we get the new cluster variables $\Omega_i=\Delta_{i-1,i+1}$, $i=2,\ldots, n-1$, and
\[
\Omega_{ij}=\Delta_{\{i-1,\ldots, i+j-2, i+j\}}\Delta_{\{i+1,\ldots, i+j-1\}} -
\Delta_{\{i,\ldots, i+j-2\}}\Delta_{\{i+1,\ldots, i+j\}}. \]

For each new quiver and the corresponding cluster variables we proceed new mutations and so on (we use Pl\"ucker
relations for calculating cluster variables as we demonstrated above).
The resulting set of cluster variables constitute the algebra $\mathcal P(n)$ being such a specialization of
$\mathcal A(n)$.  (For $n=4$ see \cite{Fomin10}, Section 2.)

Recall that  $\mathbb C[GL_n/N]$ is
the set of regular functions on $GL_n(\mathbb C)$ which are invariant under the action of $N$
by the right multiplication, that is  $\mathbb C[GL_n/N]=\mathbb C[x_{ij}]^N$.

From the next theorem follows that $\mathcal P(n)$ is a cluster algebra on $\mathbb C[GL_n/N]$.
(Recently, in \cite{GY,GY1} have shown that
the standard quantum deformation of the coordinate ring of a double Bruhat cell is a quantum cluster algebra
confirming the Berenstein-Zelevinsky conjecture \cite{BZq}. Passing to the classical limit, this shows that the
coordinate ring of a double Bruhat cell is in fact a genuine cluster algebra.)
\medskip

{\bf Theorem M0}. \begin{enumerate}
\item Any cluster variable of  $\mathcal P(n)$ belongs to  $\mathbb C[x_{ij}]^N$;
\item For every $I\subset [n]$, the flag minor $\Delta_I$ is a cluster variable of  $\mathcal P(n)$.
\end{enumerate}

For the proof of the item 2 see \cite{siena10}. For the proof of the item 1, we consider
a special seed of $\mathcal P(n)$. Namely, consider the quiver $T(n)$ corresponding to
the pseudo-line arrangement for the following reduced decomposition
(see \cite{BFZup} and \cite{siena10} for algorithms how to associate a quiver to a pseudo-line arrangement)
\[
w_0=\prod_{i=0}^{[(n-1)/2]}
s_{i+1}s_{i+2}\ldots s_{n-i-1}s_{n-i}s_{n-i-1}\ldots s_{i+2}s_{i+1}. \]

An important property of the quiver $T(n)$ is that (i) it can be reached by quiver mutations from $Q(n)$,  (ii)
the cluster variables  (in $\mathcal P(n)$)
which label the vertices of $T(n)$  are flag minors,  and (iii) any seed which is obtained by a single
mutation from $T(n)$ has properties (i) and (ii).

For $n=8$, we draw below a tiling (dual object to pseudo-line arrangement) and the corresponding quiver $T(8)$.

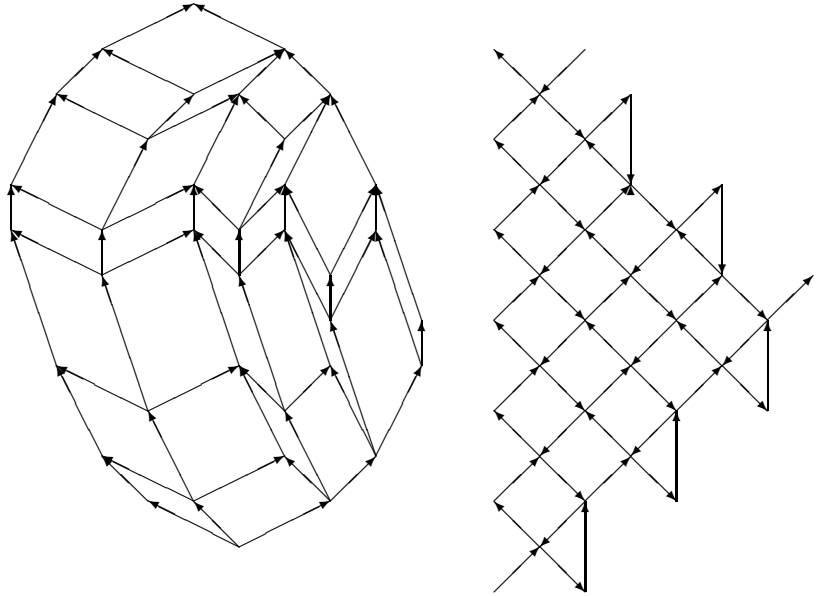
\[
\begin{array}{cc}
{
\unitlength=0.6mm \special{em:linewidth 0.4pt} \linethickness{0.4pt}
\begin{picture}(120.00,130.00)
\put(80.00,10.00){\vector(-2,1){20.00}}
\put(60.00,20.00){\vector(-1,1){10.00}}
\put(50.00,30.00){\vector(-1,2){10.00}}
\put(40.00,50.00){\vector(-1,3){10.00}}
\put(30.00,80.00){\vector(0,1){10.00}}
\put(30.00,90.00){\vector(1,2){10.00}}
\put(40.00,110.00){\vector(1,1){10.00}}
\put(50.00,120.00){\vector(2,1){20.00}}
\put(80.00,10.00){\vector(2,1){20.00}}
\put(100.00,20.00){\vector(1,1){10.00}}
\put(110.00,30.00){\vector(1,2){10.00}}
\put(120.00,50.00){\vector(0,1){10.00}}
\put(120.00,60.00){\vector(-1,3){10.00}}
\put(110.00,90.00){\vector(-1,2){10.00}}
\put(100.00,110.00){\vector(-1,1){10.00}}
\put(90.00,120.00){\vector(-2,1){20.00}}
\put(80.00,10.00){\vector(-1,1){10.00}}
\put(70.00,20.00){\vector(-1,2){10.00}}
\put(60.00,40.00){\vector(-1,3){10.00}}
\put(50.00,70.00){\vector(0,1){10.00}}
\put(50.00,80.00){\vector(1,2){10.00}}
\put(60.00,100.00){\vector(1,1){10.00}}
\put(70.00,110.00){\vector(2,1){20.00}}
\put(70.00,20.00){\vector(-2,1){20.00}}
\put(60.00,40.00){\vector(-2,1){20.00}}
\put(50.00,70.00){\vector(-2,1){20.00}}
\put(50.00,80.00){\vector(-2,1){20.00}}
\put(60.00,100.00){\vector(-2,1){20.00}}
\put(70.00,110.00){\vector(-2,1){20.00}}
\put(100.00,20.00){\vector(-1,1){10.00}}
\put(90.00,30.00){\vector(-1,2){10.00}}
\put(80.00,50.00){\vector(-1,3){10.00}}
\put(70.00,80.00){\vector(0,1){10.00}}
\put(70.00,90.00){\vector(1,2){10.00}}
\put(80.00,110.00){\vector(1,1){10.00}}
\put(70.00,20.00){\vector(2,1){20.00}}
\put(60.00,40.00){\vector(2,1){20.00}}
\put(50.00,70.00){\vector(2,1){20.00}}
\put(50.00,80.00){\vector(2,1){20.00}}
\put(60.00,100.00){\vector(2,1){20.00}}
\put(100.00,20.00){\vector(-1,2){10.00}}
\put(90.00,40.00){\vector(-1,3){10.00}}
\put(80.00,70.00){\vector(0,1){10.00}}
\put(80.00,80.00){\vector(1,2){10.00}}
\put(90.00,100.00){\vector(1,1){10.00}}
\put(110.00,30.00){\vector(-1,2){10.00}}
\put(100.00,50.00){\vector(-1,3){10.00}}
\put(90.00,80.00){\vector(0,1){10.00}}
\put(90.00,90.00){\vector(1,2){10.00}}
\put(90.00,40.00){\vector(1,1){10.00}}
\put(90.00,40.00){\vector(-1,1){10.00}}
\put(80.00,70.00){\vector(1,1){10.00}}
\put(80.00,70.00){\vector(-1,1){10.00}}
\put(80.00,80.00){\vector(1,1){10.00}}
\put(80.00,80.00){\vector(-1,1){10.00}}
\put(90.00,100.00){\vector(-1,1){10.00}}
\put(110.00,30.00){\vector(-1,3){10.00}}
\put(100.00,60.00){\vector(0,1){10.00}}
\put(100.00,70.00){\vector(1,2){10.00}}
\put(120.00,50.00){\vector(-1,3){10.00}}
\put(110.00,80.00){\vector(0,1){10.00}}
\put(100.00,60.00){\vector(1,2){10.00}}
\put(100.00,60.00){\vector(-1,2){10.00}}
\put(100.00,70.00){\vector(-1,2){10.00}}

\end{picture}

} &

{
\unitlength=0.6mm \special{em:linewidth 0.4pt} \linethickness{0.4pt}
\begin{picture}(80.00,135.00)
\put(10.00,0.00){\vector(1,1){10.00}}
\put(20.00,10.00){\vector(-1,1){10.00}}
\put(10.00,20.00){\vector(1,1){10.00}}
\put(20.00,30.00){\vector(-1,1){10.00}}
\put(10.00,40.00){\vector(1,1){10.00}}
\put(20.00,50.00){\vector(-1,1){10.00}}
\put(10.00,60.00){\vector(1,1){10.00}}
\put(20.00,70.00){\vector(-1,1){10.00}}
\put(10.00,80.00){\vector(1,1){10.00}}
\put(20.00,90.00){\vector(-1,1){10.00}}
\put(10.00,100.00){\vector(1,1){10.00}}
\put(20.00,110.00){\vector(-1,1){10.00}}
\put(30.00,120.00){\vector(-1,-1){10.00}}
\put(20.00,110.00){\vector(1,-1){10.00}}
\put(30.00,100.00){\vector(-1,-1){10.00}}
\put(20.00,90.00){\vector(1,-1){10.00}}
\put(30.00,80.00){\vector(-1,-1){10.00}}
\put(20.00,70.00){\vector(1,-1){10.00}}
\put(30.00,60.00){\vector(-1,-1){10.00}}
\put(20.00,50.00){\vector(1,-1){10.00}}
\put(30.00,40.00){\vector(-1,-1){10.00}}
\put(20.00,30.00){\vector(1,-1){10.00}}
\put(30.00,20.00){\vector(-1,-1){10.00}}
\put(20.00,10.00){\vector(1,-1){10.00}}
\put(30.00,0.00){\vector(0,1){20.00}}
\put(30.00,20.00){\vector(1,1){10.00}}
\put(40.00,30.00){\vector(-1,1){10.00}}
\put(30.00,40.00){\vector(1,1){10.00}}
\put(40.00,50.00){\vector(-1,1){10.00}}
\put(30.00,60.00){\vector(1,1){10.00}}
\put(40.00,70.00){\vector(-1,1){10.00}}
\put(30.00,80.00){\vector(1,1){10.00}}
\put(40.00,90.00){\vector(-1,1){10.00}}
\put(30.00,100.00){\vector(1,1){10.00}}
\put(40.00,110.00){\vector(0,-1){20.00}}
\put(40.00,90.00){\vector(0,0){0.00}}
\put(40.00,90.00){\vector(1,-1){10.00}}
\put(50.00,80.00){\vector(-1,-1){10.00}}
\put(40.00,70.00){\vector(1,-1){10.00}}
\put(50.00,60.00){\vector(-1,-1){10.00}}
\put(40.00,50.00){\vector(1,-1){10.00}}
\put(50.00,40.00){\vector(-1,-1){10.00}}
\put(40.00,30.00){\vector(1,-1){10.00}}
\put(50.00,20.00){\vector(0,1){20.00}}
\put(50.00,40.00){\vector(1,1){10.00}}
\put(60.00,50.00){\vector(-1,1){10.00}}
\put(50.00,60.00){\vector(1,1){10.00}}
\put(60.00,70.00){\vector(-1,1){10.00}}
\put(50.00,80.00){\vector(1,1){10.00}}
\put(60.00,90.00){\vector(0,-1){20.00}}
\put(60.00,70.00){\vector(1,-1){10.00}}
\put(70.00,60.00){\vector(-1,-1){10.00}}
\put(60.00,50.00){\vector(1,-1){10.00}}
\put(70.00,40.00){\vector(0,1){20.00}}
\put(70.00,60.00){\vector(1,1){10.00}}
\end{picture}
}

\end{array}
\] \bigskip

On the left hand side we depicted, for $n=8$, the tiling corresponding to the above defined reduced decomposition,
and on the right hand side is the quiver $T(8)$.
\medskip

For $T(n)$, the corresponding flag minors are labeled by the following sets
$\{1\}$, $\{1,2\},\ldots ,\{ 1,2\ldots n\}$,
and pairs $\{2\}, \{2,3\},\ldots, \{2,3\ldots, n-1, n\}$, $\{2, n\}, \{2,3,n\},\ldots, \{2,3\ldots, n-2,n\}$,
$\ldots$, $\{k,n-k+3,\ldots, n-1, n\}, \{k,k+1, n-k+3,\ldots, n-1, n\}\ldots , \{k,k+1,\ldots, n-k+2, n-k+3,\ldots,
n-1, n\}$,
$\{k, n-k+2,,n-k+3,\ldots, n-1, n\}, \{k,k+1, n-k+2,  n-k+3,\ldots, n-1, n\}\ldots , \{k,k+1,\ldots,n-k, n-k+2,
n-k+3,\ldots, n-1, n\}$, $k=3,\ldots, [n/2]$.
One has to assign them from the bottom to top and from the right to the left to the vertices of the tiling and
the corresponding quiver $T(n)$.
\bigskip

Let us recall  the following facts on $\mathbb C[x_{ij}]^N$ (see, for example \cite{Fulton}).
\begin{itemize}
\item  The coordinate ring $\mathbb C[x_{ij}]^N$ is a unique factorization domain.
\item  For any $I\subset [n]$, the flag minor $\Delta_I$ is irreducible polynomial in $\mathbb C[x_{ij}]$
(and hence in
$\mathbb C [x_{ij}]^N$).
\item  Irreducible flag minors $\Delta_I$ are non-equivalent among themselves. \end{itemize}

{\em Proof Theorem M0}. We prove item 1. We use the seed for the quiver $T(n)$.
Suppose there is a cluster variable $P$ which is not a polynomial in flag minors, that is
$P\not\in \mathbb C[x_{ij}]^N$.
According to Theorem CL,  $P$ is a Laurent polynomial in the variables being
the flag minors for the seed with the quiver $T(n)$. Then by assumption, this Laurent polynomial
does not belong to  $\mathbb C[x_{ij}]^N$.

Firstly, we claim that the denominator of this Laurent polynomial can not have multiplies in non-frozen variable.  In fact, suppose this not the case and the denominator has a multiplies, say $\Delta_I$,
corresponding to a non-frozen
vertex $v$  of $T(n)$. Then let us apply the mutation at this vertex. The new seed $T^v(n)$
has cluster variables being also flag minors. The latter set does not contain $\Delta_I$.
Since $P$ is a Laurent polynomial in the flag minors of the seed  $T^v(n)$, we get that the denominator
of this Laurent polynomial does not contain $\Delta_I$. This contradicts to the above facts on
$\mathbb C[x_{ij}]^N$.

Secondly, we claim that cluster variables can not have denominator with multiplies in
the frozen flag minors. We proceed by
induction on the number of mutations from the seed $T(n)$. The base of  induction is valid by
the construction of the seed $T(n)$. Consider the first case when a seed has a cluster variable with
the denominator being a monomial in frozen flag minors.
This means that for some polynomials $S$, $T$, $U$ and $V\in \mathbb C[x_{ij}]^{N}$, we have
$ \frac{S+T}{U}=\frac{V}{M}$,
where $M$ is a monomial in the frozen flag minors. Hence $M(S+T)=VU$. Because of the unique factorization property,
we have that $U$ divides $M$, and $U$ is obtained by the smaller number of mutations than $V$.
Now, we consider the first  case in the subsequence of mutations between $T(n)$ and a seed containing $U$,
such that it appears
a cluster variable which divides a monomial in frozen flag minor. That is, for some polynomials $S'$, $T'$, $U'$,
$V'\in \mathbb C[x_{ij}]^{N}$, it holds that
$
\frac{S'+T'}{U'}=V'M'$, where $M'$ is a monomial in frozen variables, and $S'$, $T'$, $U'$ do not divide frozen
variables. This contradicts that $\mathbb C[x_{ij}]^{N}$ is the unique factorization domain.
\hfill $\Box$
\medskip

We need more facts and notions to formulate our main results.

Recall, that a  partition $\lambda=(\lambda_1\ge\ldots\ge \lambda_k > 0)$ can
be identified with a Young diagram: a left-justified shape of $k$ rows of
boxes of lengths $\lambda_1,\ldots, \lambda_k$ (numbered from bottom to top for the French style).
A semi-standard Young tableau in the alphabet $[n]$ is a filling (assignment a number in $[n]$
to each box) of $\lambda$ non-decreasing  along rows and increasing along columns.
Here is an example of a semi-standard Young tableau
\begin{align*}
&\overline{|\underline 3|}\\
&{|\underline 2}\ \overline{|\underline  2}\ \overline{|\underline 2}\ \overline{|\underline  3}\
\overline{|\underline  3}\ \overline{|\underline 3 |}\\
&{|\underline 1}\ {|\underline 1}\ {|\underline 1}\ {|\underline 1}\ {|\underline 2}\ {|\underline 2}\
\overline{|\underline 2}\ \overline{|\underline 3}\ \overline{|\underline 3}\
\overline{|\underline 3}\ \overline{|\underline 3|}
\end{align*}
of the shape $(11,6,1)$.

There is a bijection between the set of semi-standard Young tableaux and integer $\mathbf D$-tight arrays
(\cite{array}).
Namely,
let us identify the space $\mathbb R^{\frac{n(n+1)}2}$  and the space of of upper
triangular $n\times n$-matrices $(a_{ij})$, $1\le j\le i\le n$.
Since we consider the French style to draw Young diagrams, we
write matrices in Descardes coordinates. Because of this, such a written matrix is
called
an {\em array}. A typical upper-triangular array $A$ is
$\begin{pmatrix} 0&0&\ldots&0 & a_{nn}\cr 0&0&\ldots& a_{n-1\,n-1} & a_{n\,n-1}\\
 \vdots&\vdots&\ldots&\vdots & \vdots\\ 0&a_{22}&\ldots&a_{n-1\, 2} & a_{n2} \\
 a_{11}&a_{21}&\ldots&a_{n-1\,1} & a_{n\,1}\end{pmatrix}$.

For a semi-standard Young tableau $T$, let us define an array $A(T)$ by the following rule: $a_{ij}(T)$ is equal
to the number of boxes in the $j$-th row filled with $i$. The array corresponding to the above tableau is
$
\begin{pmatrix} 0 & 0 & 1 \\
0 & 3 & 3 \\
4 & 3 & 4
\end{pmatrix}$.

The set of arrays corresponding the semi-standard Young tableaux is the set of integer points
in the cone of $\mathbf D$-tight
arrays.  Namely, we have (\cite{array})

{\bf Proposition}. An (integer) $n\times n$ array $A$ is an array corresponding to
a semistandard Young tableau if and only if $A$ belongs
to the cone $\mathbf D(n)$ of $\mathbf D$-{\em tight arrays}, that is, for
any $j<n$ and $i\in [n]$, it holds
\[
a(1,j+1)+\ldots +a(i,j+1)\le a(1,j)+\ldots +a(i-1,j).
\]
\medskip

A Gelfand-Tsetlin pattern is a triangular array
\[
\begin{array}{ccccccccc}
x_{nn}& & x_{n-1n}& &x_{n-2n} &&\ldots && x_{1n}\cr
&x_{n-1n-1}& &x_{n-2n-1}&&\ldots& & x_{1n-1}\cr
&&x_{n-2n-2}& &\ldots& &x_{1n-2}\cr
&&&\ddots&&&&  \cr
&&&&x_{11}&& \cr \end{array}
\]
such that $0\le x_{ij}\le x_{i-1\, j-1}\le x_{i-1\,j}$. The set of all G-T patterns is the cone $GT(n)$.

There is a linear isomorphism sending the cone $\mathbf D(n)$  to the Gelfand-Tseitlin cone $GT(n)$.
Namely, this linear isomorphism is:
$x_{11}=a_{11}$, $x_{22}=a_{22}$, $x_{1 2}=a_{11}+a_{21}$, $\ldots$, $x_{ij}=a_{ii}+\ldots +a_{ji}, \ldots,
x_{1n}=a_{11}+a_{21}+\ldots a_{n1}$.
\medskip

It is convenient to us to work with the cone $\mathbf D(n)$. We endow this cone with a tropical semi-ring structure.
Namely, the multiplication $\odot=+$ is the sum of arrays, $A\odot B:=A+B$. To define the
sum, consider  the lexicographical order and the corresponding total order on monomials
as in Introduction. Then $A\succ B$ iff $X^A$ is bigger than $X^B$ wrt the total order on monomials, where
\[
 X^A=x_{11}^{a_{11}}x_{21}^{a_{21}}\cdots x_{n1}^{a_{n1}}\cdot x_{12}^{a_{12}}x_{22}^{a_{22}}\cdots x_{n2}^{a_{n2}}
 \cdots x_{1n}^{a_{1n}}x_{2n}^{a_{2n}}\cdots x_{nn}^{a_{nn}}.
\]
Then, for a pair
of
arrays $A$, $B$ we define $\max (A,B)$ to be equal $A$ if
$A\succ B$ and be equal to $B$, otherwise. Then, we set $A\oplus
B=\max(A,B)$.

There is a distinguished  basis in $\mathbb C[x_{ij}]^N$, the tableau basis.
Elements of the {\em tableau basis} are labeled by semi-standard Young tableaux.
Namely, for a semi-standard Young tableau $T$ of shape $\lambda=(\lambda_1\ge \lambda_2\ge\ldots\ge \lambda_n)$,
we denote
$T_1,\ldots, T_{\lambda_1}$ the collections of subsets of $[n]$ being fillings of the columns of $T$.

For each $T_i$, we consider the
flag minor $\Delta_{T_i}$.
Thus, for a semi-standard Young tableau $T$, we define a monomial in flag minors as follows
\[
\Delta_{T}=\Delta_{T_1}\cdots\Delta_{T_{\lambda_1}}. \]
Let us note that the leading term of $\Delta_T$ is $X^{A(T)}$.

For the above example of semi-standard Young tableau, the corresponding monomial is \\
$\Delta_{123}(\Delta_{12})^2\Delta_{13}(\Delta_{23})^2\Delta_2(\Delta_3)^4$. \medskip

The importance of such monomials is that they form a linear basis in $\mathbb C[GL_n/N]$.\medskip

{\bf Theorem A}. (\cite{DKR, Fulton}).
The set of polynomials $\Delta_T$, while $T$ runs the set of semi-standard Young
tableaux filled from the alphabet $[n]$, form
a linear basis of $\mathbb C[x_{ij}]^N$.\medskip

{\bf Remark.} There are many other linear bases indexed by semi-standard Young tableaux. These bases
differ in rules of
forming (triangulations of $\mathbf D(n)$), for a given tableau $T$,  a monomial in flag minors (see \cite{fpsac}).
\medskip

Now we formulate main results of the paper.
\medskip

{\bf Theorem M1}. Any cluster variable $P$ in $\mathcal P(n)$ has the form
\[
P=\Delta_T+P',
\]
for some semi-standard Young tableau  $T$  and some $P'\in \mathbb C[x_{ij}]^N$
which is smaller $\Delta_T$ wrt the total
lexicographical order.
\medskip

{\bf Theorem M2}. Let $P=\Delta_T+P'$ and $Q=\Delta_T+Q'$ be cluster variables in $\mathcal P(n)$ represented
according to Theorem M1 with the same semi-standard Young tableau $T$. Then
\[
P=Q.
\]

Because of these theorems, for each seed $\mathcal S$ in  $\mathcal P(n)$, we have a collection
of semi-standard Young tableaux, $T_1(\mathcal S),\ldots, T_{(n-1)(n-2)/2}(\mathcal S)$, corresponding
to mutable variables, and
the 'frozen' one-column Young tableaux
$T_i^f:=\begin{array}{c}
\overline{|\underline{ i}|}\\ \vdots\\|\underline{ 3}|\\ |\underline{ 2}|\\|\underline{ 1}|
\end{array}$, $i=1,\ldots, n$ and
$T_{n-i}^f:=\begin{array}{c}
\overline{ |\underline{\phantom{1i} n\phantom{1i}}|}\\ |\underline{ n-1}|\\ \vdots\\ |\underline{i+1}|\\
|\underline{\phantom{1i} i\phantom{1i}}|
\end{array}$, $i=2,\ldots, n$.

Denote $C(\mathcal S)$ the cone in $\mathbf D(n)$ spanned  by the arrays corresponding
to this set of semi-standard Young tableaux for the mutable and frozen variables.
By the construction the cone $C(\mathcal S)\subset \mathbf D(n)$ is simplicial.
Moreover, the generating set of arrays of this cone is  {\em unimodular},
that is a basis in $\mathbb Z^{n(n+1)/2}$.
The unimodularity of generators of the cone $C(\mathcal S)$ follows from:
(i) due to the definition $\mathcal P(n)$, the vertices of the quiver $Q(n)$ are labeled by
the flag minors of one-column Young tableaux filled from intervals;
(ii) the the set of arrays corresponding the one-column interval tableaux is unimodular in $\mathbb Z^{n(n+1)/2}$;
(iii) cluster mutations
according to the quiver rule preserve unimodularity (see Remark after the proof of Theorem M1).

\medskip

{\bf Theorem M3}. The cones  $C(\mathcal S)$, while $\mathcal S$ runs the set of seeds in
$\mathcal P(n)$, form a simplicial fan in $\mathbf D(n)$. \medskip

{\bf Remark}. From Theorem M3 follows that cluster monomials in $\mathbb C[GL_n/N]$ are linear independent
(for skew-symmetric cluster algebras, linear independence of cluster monomials
is proven in \cite{Keller} using categorification).

For $n=3,4,5$, $\mathbb C[GL_n/N]$ is the finite cluster algebra of types $A_1$, $A_3$ and $D_6$, respectively.
The union of the cones in the fan is the whole $\mathbf D(n)$, $n=3,4,5$. Namely, we have

{\bf Corollary}. For $n=3,4,5$, the union of cones of the fan coincides with $\mathbf D(n)$.

For $n\ge 6$ the union of cones of the fan is only a part of $\mathbf D(n)$. \medskip

We end this section  with examples for $n=4,5$.
\subsection{$n=4$}

Here is the picture on the sphere (see also \cite{Fomin10})
in $\mathbb R^3$ (we consider the quotient space by the set of frozen arrays)\footnote{
Let $\phi:\mathbb R^{n(n+1)/2}\to\mathbb R^{(n-2)(n-1)/2}$ be the quotient by frozen arrays $A(T^f_i)$ and
$A(T^f_{n-i})$. Then it is easy to verify that the collection $\phi(C(\mathcal S))$,
$\mathcal S$ runs the cluster seeds, is a simplicial fan.     } illustrating the fan of cones of seed of the cluster algebra with the initial seed $Q(4)$
(such an algebra is finite $A_3$-algebra).\medskip

\unitlength=.8mm \special{em:linewidth 0.4pt}
\linethickness{0.5pt}
\begin{picture}(114.00,95.00)(-30,0)
\put(10,10){\line(1,0){100}} \put(10,10){\line(2,1){40}}
\put(10,10){\line(1,1){35}} \put(10,10){\line(2,3){50}}
\put(50,30){\line(-1,3){5}} \put(50,30){\line(1,2){10}}
\put(50,30){\line(1,0){20}} \put(50,30){\line(3,-1){60}}
\put(70,30){\line(-1,2){10}} \put(70,30){\line(1,3){5}}
\put(70,30){\line(2,-1){40}} \put(60,50){\line(-3,-1){15}}
\put(60,50){\line(3,-1){15}} \put(60,50){\line(0,1){35}}
\put(110,10){\line(-1,1){36}} \put(110,10){\line(-2,3){50}}
\qbezier[250](45,45)(52,65)(60,85)
\qbezier[250](75,45)(67,65)(60,85)
\put(50.00,30.00){\circle{2.00}} \put(70.00,30.00){\circle{2.00}}
\put(75.00,45.00){\circle{2.00}} \put(45.00,45.00){\circle{2.00}}
\put(60.00,50.00){\circle{2.00}} \put(10.00,10.00){\circle{2.00}}
\put(110.00,10.00){\circle{2.00}} \put(60.00,85.00){\circle{2.00}}
\put(67.00,24.00){\circle{2.00}}
\put(67,24){\line(1,2){4}}
\put(67,24){\line(-4,-1){58}}
\put(75.00,24.00){\makebox(0,0)[cc]{$\mathcal C_{12}$}}
\put(54.00,16.00){\makebox(0,0)[cc]{$\mathcal C_{13}$}}
\put(6.00,6.00){\makebox(0,0)[cc]{134}}
\put(114.00,6.00){\makebox(0,0)[cc]{124}}
\put(60.00,89.00){\makebox(0,0)[cc]{14}}
\put(44.00,41.00){\makebox(0,0)[cc]{13}}
\put(75.00,41.00){\makebox(0,0)[cc]{24}}
\put(50.00,27.00){\makebox(0,0)[cc]{3}}
\put(70.00,27.00){\makebox(0,0)[cc]{2}}
\put(60.00,53.00){\makebox(0,0)[cc]{23}}
\put(60.00,36.00){\makebox(0,0)[cc]{$\mathcal C_0$}}
\put(51.00,42.00){\makebox(0,0)[cc]{$\mathcal C_1$}}
\put(69.00,42.00){\makebox(0,0)[cc]{$\mathcal C_2$}}
\put(64.00,27.00){\makebox(0,0)[cc]{$\mathcal C_3$}}
\put(39.00,29.00){\makebox(0,0)[cc]{$\mathcal C_4$}}
\put(66.00,54.00){\makebox(0,0)[cc]{$\mathcal C_5$}}
\put(79.00,32.00){\makebox(0,0)[cc]{$\mathcal C_6$}}
\put(48.00,23.00){\makebox(0,0)[cc]{$\mathcal C_7$}}
\put(54.00,54.00){\makebox(0,0)[cc]{$\mathcal C_8$}}
\put(81.00,45.00){\makebox(0,0)[cc]{$\mathcal C_9$}}
\put(39.00,45.00){\makebox(0,0)[cc]{$\mathcal C_{10}$}}
\put(95.00,48.00){\makebox(0,0)[cc]{$\mathcal C_{11}$}}
\end{picture}

There are 14 cones (the cone $\mathcal C_{11}$ is the triangle $14, 134, 124$ on the 2-sphere), 9 variables, 8 of them are of the form $A_I$, where $I$ is a non-frozen subset of $[4]$ ($2^4-1=15-7=8$) and
one more variable  corresponding to the vertex, indexed by the Young tableau $\begin{array}{cc}
\overline{|\underline 4 |}&\\\overline{|\underline 2 |}&\\\overline{|\underline1}&
\!\!\!\!\!\!\overline{|\underline 3|}\end{array}$. This Young tableau is a vertex in
cones $\mathcal C_3$, $\mathcal C_7$, $\mathcal C_{12}$, $\mathcal C_{13}$.

In the quotient space by the frozen arrays,  the images of the arrays  2, 3,  24, 14,  23, 13, 134, and
 124 are related by the following relations:

$\overline{3}+\overline{14}=0$ (because there holds 3+14=1+34 in  $\mathbf D(n)$);
$\overline{2}+\overline{134}=0$ (2+134=1+234);  $\overline{23}+\overline{124}=0$ (23+124=12+234);
$\overline{2}+\overline{14}=\overline{24}$ (2+14=1+24);
$\overline{23}+\overline{134}=\overline{13}$ (23+134=13+234).

For unique non-extreme ray cluster variable we have
      $\overline{\begin{array}{cc}
4&\\2&\\1&3\end{array}}=\overline{124}+\overline{3}$.

      Let the vectors  $\overline 3$, $\overline {23}$ and
$\overline 2$ be a basis of the 3-dim quotient space. Then these cluster vectors  are of
the form\medskip

      $\overline{3}=(1,0,0)$, $\overline{14}=(-1,0,0)$ ,
      $\overline{23}=(0,1,0)$, $\overline{124}=(0,-1,0)$,
      $\overline{2}=(0,0,1)$, $\overline{ 134}=(0,0,-1)$,
     $\overline{ 24}=(-1,0,1)$, $\overline{13}=(0,1,-1)$, and
     $\overline{\begin{array}{cc}
4&\\2&\\1&3\end{array}}=(1,-1,0)$.\medskip

This fan is the normal fan to the associahedron

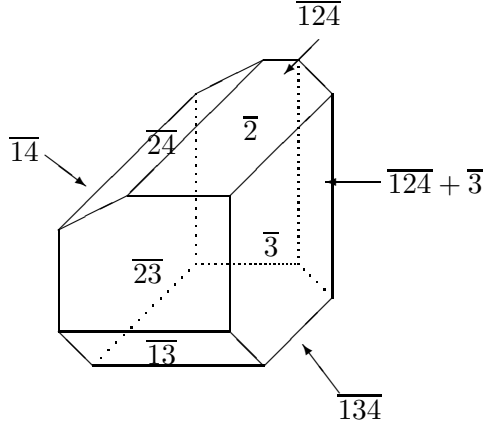
\unitlength=.9mm \special{em:linewidth 0.4pt}
\linethickness{0.4pt}
\begin{picture}(98.00,67.00)

\put(55,10){\line(-1,1){5}} \put(50,15){\line(0,1){15}}
\put(55,10){\line(1,0){25}} \put(50,15){\line(1,0){25}}
\put(50,30){\line(2,1){10}}
\put(80,10){\line(-1,1){5}}  \put(80,10){\line(1,1){10}}
\put(75,35){\line(-1,0){15}} \put(75,35){\line(0,-1){20}}
\put(75,35){\line(1,1){15}}
 \put(55,10){\line(-1,1){}}
\put(80,55){\line(-2,-1){10}}
\put(80,55){\line(1,0){5}}
\put(90,50){\line(0,-1){30}}
\put(90,50){\line(-1,1){5}}

 \put(50,30){\line(1,1){20}}  \put(80,55){\line(-1,-1){20}}

\put(105.00,37.00){\makebox(0,0)[cc]{$\overline{124}+\overline 3$}}
\put(97.00,37.00){\vector(-1,0){8.00}}

\bezier{5}(90.00,20.00)(87.00,23.00)(85.00,25.00)
\bezier{30}(85.00,55.00)(85.00,40.00)(85.00,25.00)

\bezier{20}(70.00,50.00)(70.00,39.00)(70.00,25.00)
\bezier{20}(70.00,25.00)(80.00,25.00)(85.00,25.00)
\bezier{14}(70.00,25.00)(63.00,18.00)(55.00,10.00)
\put(63.00,24.00){\makebox(0,0)[cc]{$\overline{23}$}}
\put(81.00,28.00){\makebox(0,0)[cc]{$\overline{3}$}}
\put(78.00,45.00){\makebox(0,0)[cc]{$\overline{2}$}}
\put(65.00,43.00){\makebox(0,0)[cc]{$\overline{24}$}}
\put(65.00,12.00){\makebox(0,0)[cc]{$\overline{13}$}}
\put(88.00,62.00){\makebox(0,0)[cc]{$\overline{124}$}}
\put(94.00,4.00){\makebox(0,0)[cc]{$\overline{134}$}}
\put(45.00,42.00){\makebox(0,0)[cc]{$\overline{14}$}}
\put(88.00,58.00){\vector(-1,-1){5.00}}
\put(91.00,6.00){\vector(-3,4){5.33}}
\put(48.00,41.00){\vector(4,-3){6.00}}
\end{picture}


\medskip

\subsection{n=5}
For $n=5$ the fan has 672 simplicial cones and 36  one-dimensional rays corresponding to mutable variables
and 9 rays corresponding to frozen.
The quotient of this fan by frozen rays is the normal fan to the $D_6$-associahedron.
This is 6-dimensional  polytope  which  has
672 vertices and 36 facets.

The mutable variables have the leading term $\Delta_T$ (Theorem M1) of the following form:
either $T$ is an one column semistandard Young tableaux whose filling is a subset $I\subset [5]$, such that
$I$ is not an interval containing either $1$ or $n$, (in such a case, cluster variables coincide
with the leading terms) or $T$ is one of the following  set of two column tableaux:

$\begin{array}{cc} 4&\\2&\\1&3\end{array}$, $\begin{array}{cc}
5&\\2&\\1&3\end{array}$, $\begin{array}{cc}
5&\\2&\\1&4\end{array}$, $\begin{array}{cc}
5&\\3&\\2&4\end{array}$, $\begin{array}{cc}
5&\\3&\\1&4\end{array}$, $\begin{array}{cc}
4&\\2&5\\1&3\end{array}$, $\begin{array}{cc}
5&\\3&4\\1&2\end{array}$, $\begin{array}{cc}
5&\\3&\\2&\\1&4\end{array}$, $\begin{array}{cc}
5&\\3&\\2&4\\1&3\end{array}$, $\begin{array}{cc}
5&\\3&\\2&4\\1&2\end{array}$, $\begin{array}{cc}
5&\\3&\\2&4\\1&1\end{array}$,  $\begin{array}{cc}
5&\\4&\\2&\\1&3\end{array}$, $\begin{array}{cc}
5&\\4&\\2&4\\1&3\end{array}$, $\begin{array}{cc}
5&\\4&\\2&5\\1&3\end{array}$.
\medskip

\section{Proof of Theorems M1-M3}

{\em Proof of Theorem M1}.
Because of Theorem M0 (item 1), each cluster variable is a polynomial in flag minors.
So, we have to prove that each such cluster polynomial $P$ has the coefficient
$1$ with the leading monomial wrt the lexicographic
order.
Proceed by induction on the number of mutations from the initial seed with the quiver $Q(n)$. For
one mutations we gave explicit formulae for $\Omega_{ij}$ and $\Omega_i$ above, and we are done.

Suppose we are done for some number of mutations, and consider one more mutation.
Consider the quiver mutation rule
\[
x^{new}_v\cdot x_v=\prod_{v'\in In(v)}x_{v'}^{w(v'v)}+\prod_{v''\in
Out(v)}x_{v''}^{w(v,v'')}.
\] By induction, for each $v$, $v'\in In(v)$, and $v''\in Out (v)$,
we have $x_v=\Delta_{T(v)}+P'_v$, $x_{v'}=\Delta_{T(v')}+P'_{v'}$ and $x_{v''}=\Delta_{T(v'')}+P'_{v''}$.
By Theorem M0, $x_v^{new}=a\Delta_{T}+P'$. Then, because the product of leading terms on the left hand side equals
the product of the leading terms on the right hand side, we have
\[
aX^{A(T)}X^{A(T(v))}=X^{\max(\sum_{v'\in I(v)}w(v',v)A(T(v')), \sum_{v''\in Out(v)}w(v,v'')A(T(v'')))}, \]
where the operation $\max$ for the cone $\mathbf D(n)$ is defined above.

Because of this, we have $a=1$.\hfill $\Box$

{\bf Remark}.
Let us note that from the proof of this theorem, for any seed of $\mathcal P(n)$ and any vertex $v$,
we have
 \begin{equation}\label{tropical1}
   \max(\sum_{v'\in I(v)}w(v',v)A(T(v')), \sum_{v''\in Out(v)}w(v,v'')A(T(v'')))-A(T(v)\in \mathbf D(n).
 \end{equation}

From this relation follows that if, for a seed $\mathcal S$, the
set of arrays $\{A(T_v)\}$ (vectors in $\mathbb Z^{n(n+1)/2}$),
$v$ runs the vertices of the quiver for $\mathcal S$,
is unimodular set of vectors in $\mathbb Z^{n(n+1)/2}$. Then, for any
seed, which is obtained by a mutation, and hence, for all seeds in $\mathcal P(n)$, the unimodularity
is valid. In fact, due to (\ref{tropical1}),
the transformation matrix for mutations of cluster cones is unimodular.
For the initial quiver $Q(n)$, the interval tableaux, the labels of the vertices of $Q(n)$, form a unimodular set.
Therefore every cone $C(\mathcal S)$ is unimodular. \bigskip

{\em Proof of Theorem M2}. Suppose  there are two seeds $\mathcal S$ and ${\mathcal S}'$
and cluster variables $w$ and $w'$ in these seeds such that $w=\Delta_T+P$ and $w'=\Delta_T+P'$, where
$P$ and $P'$ are polynomials in the flag minors which are smaller (wrt the lexicographical order) $\Delta_T$.

Suppose $w\neq w'$. By Theorem CL,  $w$ is a Laurent polynomial of the cluster variables of
the seed ${\mathcal S}'$.
We claim that this  Laurent polynomial has the form
\begin{equation}\label{sum1}
w=w'+p'( S'\setminus w'),
\end{equation}
where $p'( S'\setminus w')$ is a Laurent polynomial in variables of the seed $\mathcal S'$ except $w'$.
In fact, the leading term of $w$ is $X^{A(T)}$. Then there exists a  Laurent monomial in variables $S'$ of the above
Laurent polynomial, such that its leading term (that is the product of leading term of this monomial) is equal to
$X^{A(T)}$. But the leading term of $w'$ is also $X^{A(T)}$. Hence there is a linear dependence between
arrays for tableaux corresponding to the variables of $\mathcal S'$. Due to Remark after proof Theorem M1, this
is no the case. Thus (\ref{sum1}) holds true.

By the same line of arguments, it holds that
\[
w'=w+p(S\setminus w),
\]
where $p(S\setminus w)$ is a Laurent polynomial in variables of the seed $\mathcal S\setminus w$.

Therefore, we have
\[
p(S\setminus w)+p'(S'\setminus w')=0.
\]
Because of the  positivity of coefficients of Laurent polynomials (Theorem CP) , the latter equality is possible if and only if
\[
p=p'=0.
\]
In fact, the numerator of $p(S\setminus w)+p'(S'\setminus w')$
is a Laurent polynomial in variables of the initial seed (with the quiver
$Q(n)$) with positive coefficients. Since, for a totally
positive matrix,  all interval minors are positive we get $p(S\setminus w)+p'(S'\setminus w')>0$.
Thus, $p=p'=0$, and hence $w=w'$.
\hfill  $\Box$
 \medskip

{\em Proof of Theorem M3}.

Suppose there are two seeds, $\mathcal S$ and $\mathcal S'$,
such that the cones $C:=C(\mathcal S)$ and $C':=C(\mathcal S')$  intersect   
such that it holds $int(C)\cap int(C')\neq\emptyset$ ($int(C')$ denotes interior of a cone).
Then, because the cones are unimodular,
this intersection contains an integer array, say, corresponding a semi-standard Young tableau $T$. Then
there is a cluster monomial $y$  in the variables of the seed $\mathcal S$ such that $y=x^{A(T)}+$
lexicographically smaller terms,
and there is a cluster monomial $y'$  in the variables of the seed $\mathcal S'$ such that
$y'=x^{A(T)}+$ lexicographically smaller terms.

Let us show that $y=y'$.
By Theorem CL,
the cluster monomial  $y$ is the product of Laurent polynomials  in the variables of the seed $\mathcal S'$.
The same reasoning as in the proof Theorem M2 shows that $y=y'+R(S')$. Similarly, we get $y'=y+R'(S)$,
where $R$ and $R'$ are Laurent polynomials with positive coefficients (Theorem CP). Thus, $R(C')+R(C)=0$ is
possible iff $R(C')=R(C)=0$. Hence $y=y'$.

Now, recall that $y=\prod_{v}s_v^{a_v}$ is the monomial in the cluster variables of $S$.
By Theorem CL, $y$ is the product of
the Laurent polynomials in cluster variables of $S'$,
since each $s_v$ is a Laurent polynomial in $S'$. But $y=y'$ and $y'$ is a monomial in $S'$. Hence each cluster
variable $s_v$ is a Laurent monomial in $S'$. And, similarly, each cluster variable in $S'$ is a Laurent monomial
in $S$.

Since the cones $C:=C(\mathcal S)$ and $C':=C(\mathcal S')$ have a common interior point, there is a facet
of one of this cones which has a common interior point with another cone. Let this facet be a facet of
$C:=C(\mathcal S)$. Note that the vertices corresponding to the frozen arrays belong to this facet. Let
$\mathcal S''$ be a seed obtained through the mutation of $\mathcal S$ in the vertex, corresponding to the
ray of $C:=C(\mathcal S)$ which does not belong to this facet.
Then $C'':=C(\mathcal S'')$ and $C':=C(\mathcal S')$ have a common interior point (for example, in neighbor
of the facet).   Hence, by the same reasoning as above, the new cluster variable in $\mathcal S''$
is a Laurent monomial in
cluster variable in $S'$, but this is not the case, because all cluster variables in $\mathcal S$
are  Laurent monomials in $S'$.
\hfill $\Box$

\section{Lusztig basis}

Let $X\in Mat_q(n)$ be a $n\times n$ quantum $q^2$-matrix. That is, the following relations
hold:\medskip

$x_{il}x_{ik} = q^2x_{ik}x_{il}$ ~~$\forall\,i,~ \forall\,k < l$;
$x_{jk}x_{ik} = q^2x_{ik}x_{jk}$ ~~$\forall\,i < j, ~\forall\,k$;
$x_{jk}x_{il} = x_{il}x_{jk}$ ~~$\forall\,i < j,~\forall\, k < l$;
$x_{jl}x_{ik} = x_{ik}x_{jl} + (q^2 - q^{-2})x_{il}x_{jk}$ ~~$\forall\,i < j,
~\forall\,k < l$. \medskip

Consider monomials $X^A$, $A\in \mathbf D(n)$, recall that the product is taken
w.r.t. the lexicographical order $\succ$. Define
$l(A):={\sum_i(\sum_{j>k}a_{ij}a_{ik}+\sum_{j>k}a_{ji}a_{ki})}$ and   set
$X(A):=q^{-l(A)}X^A$.

Lusztig in (\cite{L}) proved that there exists a unique basis in $\mathbb C_{q}[GL_n/N]$ (the dual canonical basis) $\mathbf
B^*(n)=\{b_A, \, A\in \mathbf D(n)\}$ determined by the following conditions:
\begin{itemize}
\item $\overline{b_A(X)}=b_A(X)$, where $\overline {qX_{ij}X_{lk}}=
q^{-1}X_{lk}X_{ij}$;
\item $b_A(X)=X(A)+\sum_{B\prec A}c_A(B)X(B)$, where $B\prec A$ denotes the
above considered lexicographical order on the set of the arrays,
and the summation is taken over the matrices $B$ with the same column and
row sums as those in $A$, and $c_A(B)\in q\mathbb Z [q]$.
\end{itemize}

Because of this definition, for $q=1$, $b_A(X)=\Delta_{T(A)}+b'(X)$, where $T(A)$ is the semi-standard
Young tableau corresponding to $A$, and $\Delta_{T(A)}\succ b'$.
Thus, the Lusztig basis is another basis of $\mathbb C[GL_n/N]$ labeled by the semi-standard Young tableaux
the (integer points of the cone $\mathbf D(n)$).

Let us say that an element $b\in \mathbf{B}^*(n)$ is {\em real} if, for any $r$, $b^r\in
\mathbf{B}^*(n)$ up to a power of $q$. In such a case, $b^r_A=q^{rl(A)}b_{rA}$.
\medskip

We call a semi-standard Young tableau {\em real} if it labels the real element of  $\mathbf {B}^*(n)$.\medskip

{\bf Conjecture}.  The labels of the cluster monomials in $\mathbb C[GL_n/N]$ constitute the set of real
semi-standard Young tableaux.

\end{document}